\documentclass{amsart}
\usepackage{amsmath,amsthm,amsopn,amstext,amscd,amsfonts,amssymb,mathrsfs,mathtools, esint}
\usepackage{dsfont}
\usepackage{comment}
\usepackage{graphicx}

\theoremstyle{theorem}

\theoremstyle{definition}

\begin{document}
\title[Remark on ``When Are All the Zeros of a Polynomial Real and Distinct?'']{Remark on ``When Are All the Zeros of a Polynomial Real and Distinct?''}
\author{K. Castillo}
\address{CMUC, Department of Mathematics, University of Coimbra, 3001-501 Coimbra, Portugal}
\email{kenier@mat.uc.pt}

\maketitle

In  \cite{C20}, Chamberland proved the following theorem: 
\begin{changemargin}{0.8cm}{0.8cm}
{\em Let $P$ be a polynomial of degree $n\geq 2$ with real coefficients. Then the zeros of $P$ are real and distinct if and only if
\begin{align}\label{eq1}
\Big(\big(P^{(n-j-1)}(x)\big)'\Big)^2-P^{(n-j-1)}(x)\big(P^{(n-j-1)}(x)\big)''>0
\end{align}
for all $x\in \mathbb{R}$ and $j\in\{1,\dots, n-1\}$.
}
\end{changemargin}

The purpose of this note is to point out that this result is implicitly contained in the elementary lore of the theory of orthogonal polynomials on the real line (OPRL).
 
$(\Rightarrow)$ By Geronimus-Wendroff's theorem \cite[Theorem 4.4.19]{S15}, there are probability measures on $\mathbb{R}$, $\mu_j$, with finite moments so that $P_{j+1}=P^{(n-j-1)}$ and $P_j=P'_{j+1}$ are among the OPRL for $\mu_j$. Since $P_{j+1}$ and $P_j$ have leading coefficients of the same sign, then \cite[(3.3.6)]{S59}
\begin{align}\label{eq2}
P'_{j+1}(x)P_j(x)-P'_j(x)P_{j+1}(x)>0
\end{align}
for all $x\in \mathbb{R}$ and $j\in\{1,\dots, n-1\}$.

$(\Leftarrow)$ Set $P_{j}=P^{(n-j)}$ for all $j\in \{1,\dots, n\}$. Without loss of generality we can assume $P$ monic. Suppose that the zeros of $P_j$ are real and distinct. From \eqref{eq2} we see that $P_{j+1}$ has at least one zero between two zeros of $P_{j}$, and so $P_{j+1}$ has at least $j-1$ real and distinct zeros. Suppose that the other two zeros of $P_{j+1}$ are not real, and therefore they appear as a complex conjugate pair.  Let $a$ be the largest zero of $P_j$. Clearly, $P_{j+1}(x)>0$ for all $x>a$. However, by \eqref{eq2}, we have $P_{j+1}(a)<0$, which leads to a contradiction. From this we conclude that if $P_{j}$ has real and distinct zeros, then the same holds for $P_{j+1}$.  Since $P_1=P^{(n-1)}$ has a single real zero, we infer successively that $P_2, \dots, P_n=P$ have real and distinct zeros.

We emphasize for the reader's convenience that any polynomial with real and distinct zeros is an element of a sequence of OPRL, and so any sentence starting with {\em ``the zeros of a polynomial are real and distinct if and only if"} is virtually talking about OPRL.

\end{document}